\documentclass[12pt]{amsart}
\usepackage{amscd}
\def\NZQ{\Bbb}               % the font for N,Z,Q,R,C

\def\ZZ{{\NZQ Z}}

\def\PP{{\NZQ P}}

\def\frk{\frak}               % font for "Fraktur"

\def\mm{{\frk m}}

\def\Phi{{\frk n}}
\def\Phi{{\frk N}}
\def\opn#1#2{\def#1{\operatorname{#2}}} % to make operators
\opn\chara{char} \opn\length{\ell} \opn\pd{pd} \opn\rk{rk}
\opn\projdim{proj\,dim} \opn\injdim{inj\,dim} \opn\rank{rank}
\opn\depth{depth} \opn\and{and} \opn\grade{grade}
\opn\height{height} \opn\embdim{emb\,dim} \opn\codim{codim}

\opn\Tr{Tr} \opn\bigrank{big\,rank}
\opn\superheight{superheight}\opn\lcm{lcm}
\opn\trdeg{tr\,deg}%
\opn\reg{reg} \opn\lreg{lreg} \opn\ini{in}
\opn\div{div} \opn\Div{Div} \opn\cl{cl} \opn\Cl{Cl}
\opn\Spec{Spec} \opn\Supp{Supp} \opn\supp{supp} \opn\Sing{Sing}
\opn\Ass{Ass} \opn\Min{Min}
\opn\Ann{Ann} \opn\Rad{Rad} \opn\Soc{Soc}
\opn\Im{Im}
 \opn\Ker{Ker} \opn\Coker{Coker} \opn\Am{Am}
\opn\Hom{Hom} \opn\Tor{Tor} \opn\Ext{Ext} \opn\End{End}
\opn\Aut{Aut} \opn\id{id}

\opn\nat{nat}
\opn\pff{pf}%   \pf exists already
\opn\Pf{Pf} \opn\GL{GL} \opn\SL{SL} \opn\mod{mod} \opn\ord{ord}
\opn\cl{cl} \opn\conv{conv} \opn\ext{ext} \opn\rad{rad}
\opn\red{red}
\opn\aff{aff} \opn\con{conv} \opn\relint{relint} \opn\st{st}
\opn\lk{lk} \opn\cn{cn} \opn\core{core} \opn\vol{vol}
\opn\link{link} \opn\star{star}
\opn\gr{gr}

\def\pot#1#2{#1[\kern-0.28ex[#2]\kern-0.28ex]}

\opn\dirlim{\underrightarrow{\lim}}
\opn\inivlim{\underleftarrow{\lim}}
\let\iso=\cong

\let\Dirsum=\bigoplus

\let\to=\rightarrow

\def\Implies{\ifmmode\Longrightarrow \else
     \unskip${}\Longrightarrow{}$\ignorespaces\fi}
\def\implies{\ifmmode\Rightarrow \else
     \unskip${}\Rightarrow{}$\ignorespaces\fi}
\def\iff{\ifmmode\Longleftrightarrow \else
     \unskip${}\Longleftrightarrow{}$\ignorespaces\fi}

\let\:=\colon
\newtheorem{Theorem}{Theorem}[section]

\newtheorem{Corollary}[Theorem]{Corollary}
\newtheorem{Proposition}[Theorem]{Proposition}

\let\epsilon\varepsilon
\let\phi=\varphi
\let\kappa=\varkappa
\def\qed{\ifhmode\textqed\fi
   \ifmmode\ifinner\quad\qedsymbol\else\dispqed\fi\fi}
\def\textqed{\unskip\nobreak\penalty50
    \hskip2em\hbox{}\nobreak\hfil\qedsymbol
    \parfillskip=0pt \finalhyphendemerits=0}
\def\dispqed{\rlap{\qquad\qedsymbol}}

\opn\dis{dis}
\def\pnt{{\raise0.5mm\hbox{\large\bf.}}}

\begin{document}

\title{Failure of  tameness for local cohomology}

\author{Steven Dale Cutkosky  and J\"urgen Herzog}\thanks{The first author
was partially supported by NSF}

\address{Steven Dale Cutkosky, Department of Mathematics,
University of Missouri, Columbia, MO 65211,
USA}\email{CutkoskyS@missouri.edu}

\address{J\"urgen Herzog, Fachbereich Mathematik und
Informatik, Universit\"at Duisburg-Essen, Campus Essen, 45117
Essen, Germany} \email{juergen.herzog@uni-essen.de}

\maketitle

\begin{abstract}
We give an example that shows that not all local cohomology
modules are tame in the sense of Brodmann and Hellus.
\end{abstract}

\section*{Introduction}
\noindent In their paper ``Cohomological patterns of coherent
sheaves over projective schemes", Brodmann and Hellus \cite{BrHe}
raised the following tameness problem: let $R=\Dirsum_{n\geq
0}R_n$ be a positively graded Noetherian ring such that $R_0$ is
semilocal, and let $M$ be a finitely generated graded $R$-module.
Denote by $J$ the graded ideal $\Dirsum_{n>0}R_n$. Is true that
all the local cohomology modules $H_J^i(M)$ are tame? The authors
call a graded $R$-module $T$ tame, if there exists an integer
$n_0$ such that $T_n=0$ for all $n\leq n_0$, or else $T_n\neq 0$
for all $n\leq n_0$.

The tameness problem has been answered in the affirmative in many
cases, in particular, if $\dim R_0\leq 2$. We refer to the article
\cite{Br} of Brodmann for a survey on this problem. In this paper
we present an example that shows that the tameness problem does
not always have a positive answer. In our example, $R_0$ is a
$3$-dimensional normal local ring with isolated singularity.

\section{A bigraded ring with periodic local cohomology}
Suppose that ${\mathcal A}$ is a very ample line bundle on a
projective space $\PP^n$, and ${\mathcal F}$ is a coherent sheaf
on $\PP^n$. ${\mathcal F}$ is $m$-regular with respect to
${\mathcal A}$ if $H^i(\PP^n,{\mathcal F}\otimes{\mathcal
A}^{\otimes (m-i)})=0$ for all $i>0$.

If ${\mathcal F}$ is $m$-regular with respect to ${\mathcal A}$,
then $H^0(\PP^n,{\mathcal F}\otimes{\mathcal A}^{\otimes k})$ is
spanned by
$$
H^0(\PP^n,{\mathcal F}\otimes{\mathcal A}^{\otimes (k-1)})\otimes
H^0(\PP^n,{\mathcal A})
$$
if $k>m$ (Lecture 14 \cite{Mu}).

\begin{Theorem}
\label{main} Suppose that $\bf k$ is an algebraically closed
field. Then there exists a bigraded domain
$$
R=\sum_{m,n\ge 0}R_{m,n}t^mu^n
$$
with the following properties:
\begin{enumerate}
\item $R$ is of finite type over $R_{0,0}={\bf k}$, and is
generated in degree 1 over $R_{0,0}$ (with respect to the grading
$d(f)=m+n$ for $f\in R_{m,n}t^mu^n$). \item $R$ has dimension 4,
is normal, and the singular locus of $\Spec(R)$ is the bigraded
maximal ideal $\sum_{m+n>0}R_{m,n}t^mu^n$. \item For $n\ge 0$, let
$$
R_n=\sum_{m\ge 0}R_{m,n}t^m.
$$
$R_0$ is a normal, 3 dimensional graded ring, the $R_n$ are graded
$R_0$ modules, and $R=\sum_{n \ge 0}R_nu^n$ is generated by $R_1$
as an $R_0$ algebra. \item  Let $I=\sum_{m>0}R_{m,0}t^m$ be the
graded maximal ideal of $R_0$. The singular locus of $\Spec(R_0)$
is $I$. \item For $n\ge 0$, we have
$$
H^2_I(R_n)=\left\{\begin{array}{ll}
{\bf k}^2&\text{ if $n$ is even}\\
0&\text{ if $n$ is odd.}
\end{array}\right.
$$
\end{enumerate}
\end{Theorem}

\begin{proof} Let $E$ be an elliptic curve over ${\bf k}$, and $\overline D$ be a degree 0 divisor on $E$ such that
$\overline D$ has order 2 in the Jacobian of $E$ ($\overline
D\not\sim 0$ and $2\overline D\sim 0$). Let $p\in E$ be a (closed)
point. Let $S=E\times_{\bf k} E$, with projections
$\pi_1:S\rightarrow E$ and $\pi_2:S\rightarrow E$. Let
$H=\pi_1^*(p)+\pi_2^*(p)$ and $D=\pi_2^*(\overline D)$ be divisors
on $S$. $H$ and $D+H$ are ample on $S$ by the Nakai criterion
(Theorem 5.1 \cite{Ha}).

Suppose that $m,n\in {\bf Z}$.
$$
\begin{array}{ll}
H^1(S,{\mathcal O}_S(mH+nD))\cong &H^0(E,{\mathcal O}_E(mp))\otimes_{\bf k} H^1(E,{\mathcal O}_E(mp+n\overline D)) \\
&\oplus H^1(E,{\mathcal O}_E(mp))\otimes_{\bf k} H^0(E,{\mathcal
O}_E(mp+n\overline D))
\end{array}
$$
by the K\"unneth formula (IV of Lecture 11 \cite{Mu}). If $m<0$, then
$$
H^0(E,{\mathcal O}_E(mp))=0 \text{ and }H^0(E,{\mathcal
O}_E(mp+n\overline D))=0.
$$
If $m>0$, then by Serre duality,
$$
H^1(E,{\mathcal O}_E(mp))\cong H^0(E,{\mathcal O}_E(-mp))=0
$$
and
$$
H^1(E,{\mathcal O}_E(mp+n\overline D))\cong H^0(E,{\mathcal
O}_E(-mp-n\overline D))=0.
$$
If $m=0$, we have
$$
H^1(S,{\mathcal O}_S(nD))\cong H^1(E,{\mathcal O}_E(n\overline
D))\oplus H^0(E,{\mathcal O}_E(n\overline D)).
$$
By the Riemann-Roch theorem on $E$, we have for $m,n\in {\bf Z}$,
\begin{equation}\label{eq1}
h^1(S,{\mathcal O}_S(mH+nD))=\left\{\begin{array}{ll}
2&\text{ if $m=0$ and $n$ is even}\\
0&\text{ otherwise.}
\end{array}
\right.
\end{equation}

Let ${\mathcal E}={\mathcal O}_S\oplus {\mathcal O}_S(D)$. Let
$X={\bf P}({\mathcal E})$ be the  projective space bundle
(\cite[Section II.7]{Ha}) with projection $\pi:X\rightarrow S$ and
associated line bundle ${\mathcal O}_X(1)$. Since $H$ is ample on
$S$, there exists a number $\tau$ such that $\pi^*{\mathcal
O}_S(nH)\otimes_{{\mathcal O}_X}{\mathcal O}_X(1)$ is very ample
on $X$ for all $n\ge \tau$ ( \cite[Proposition II.7.10 and
Exercise II.5.12(b)]{Ha}).

There exists $l\ge \tau$ such that $lH$ and $D+lH\sim l(D+H)$ are
very ample.

There exists an odd number $r_2>0$ such that $H^i(S,{\mathcal
O}_S((r_2-i)(D+alH))=0$ for all $a\ge 1$ and $i>0$. Thus for all
$a\ge 1$, ${\mathcal O}_S$ is $r_2$ regular for $D+alH$. It
follows that $H^0(S,{\mathcal O}_S(mr_2(D+alH))$ is spanned by
$H^0(S,{\mathcal O}_S(r_2(D+alH)))^{\otimes m}$ for all $a\ge 1$
and $m\ge 1$.

There exists $r_1>0$ such that ${\mathcal O}_S$ and ${\mathcal
O}_S(D)$ are $r_1$ regular for $lH$. Thus ${\mathcal O}_S(nD)$ is
$r_1$ regular for $lH$ for all $n\in \ZZ$ (since $2D\sim 0$).

Choose $a>r_1$. For $m,n\ge 0$, $H^0(S,{\mathcal
O}_S(mr_2(D+alH)+nr_2alH))$ is spanned by $H^0(S,{\mathcal
O}_S(mr_2(D+alH))\otimes H^0(S,{\mathcal O}_S(r_2alH))^{\otimes
n}$ since ${\mathcal O}_S(mr_2 D)$ is $r_1$ regular for $lH$.

Thus $H^0(S,{\mathcal O}_S(mr_2(D+alH)+nr_2alH))$ is spanned by
$$
H^0(S,{\mathcal O}_S(r_2(D+alH))^{\otimes m}\otimes
H^0(S,{\mathcal O}_S(r_2alH))^{\otimes n}.
$$

 Let
$$
R_{m,n}=H^0(S,{\mathcal O}_S(mr_2laH+nr_2(D+alH))).
$$
Let
$$
R_n=\sum_{m\ge 0} R_{m,n}t^m,
$$
$$
R=\sum_{m,n\ge 0}R_{m,n}t^mu^n=\sum_{n\ge 0}R_nu^n,
$$
where $t,u$ are indeterminates.
  $R$ is normal (by \cite[Theorem
4.2]{Z}). By our construction, $R$ is generated as an
 $R_{0,0}={\bf k}$ algebra by $R_{1,0}t$ and $R_{0,1}u$. We deduce
 that
$R$ is standard graded as an $R_{0}$ algebra, $R_0$ is normal and
standard graded as an $R_{0,0}={\bf k}$ algebra. We have an
isomorphism $S\cong \text{Proj}(R_0)$, with associated line bundle
${\mathcal O}_S(1)\cong {\mathcal O}_S(r_2alH)$. Since $S$ is
nonsingular, the singular locus of $\Spec(R_0)$ is the graded
maximal ideal $I=\sum_{m>0}R_{m,0}t^m$ of $R_0$.  The sheafication
of the module $R_n$ on $S$ is
$$
\tilde R_n={\mathcal O}_S(nr_2(D+alH)),
$$
by Exercise II.5.9 (c) \cite{Ha}. For $n\ge 0$, by (\ref{eq1}),
and the exact sequences relating local cohomology and global
cohomology (\cite[A.4.1 ]{E}), we have
$$
H^2_I(R_n)\cong\oplus_{m\in{\bf Z}}H^1(S,{\mathcal O}_S(
mr_2alH+nr_2(D+alH)))\cong H^1(S,{\mathcal O}_S(nr_2D)).
$$
Thus
\begin{equation}\label{eq2}
H^2_I(R_n)\cong\left\{\begin{array}{ll}
{\bf k}^2&\text{ if $n$ is even,}\\
0&\text{ if $n$ is odd.}
\end{array}\right.
\end{equation}

Let ${\mathcal L}=\pi^*{\mathcal O}_S(r_2alH)\otimes_{{\mathcal
O}_X}{\mathcal O}_X(1)$. By our choice of $l$,
 ${\mathcal L}$
 is very ample on $X$. For $r\ge 0$, we have (by \cite[Proposition II.7.11 ]{Ha})
$$
\begin{array}{ll}
H^0(X,{\mathcal L}^r)&\cong H^0(S,\text{Sym}^r({\mathcal E})\otimes_{{\mathcal O}_S}{\mathcal O}_S(rr_2alH))\\
&\cong \oplus_{n=0}^rH^0(S,{\mathcal O}_S(rr_2alH+nD))\\
&\cong\oplus_{n=0}^rH^0(S,{\mathcal O}_S(rr_2alH+nr_2D))\text{ (since $r_2$ is odd and $D$ has order 2)}\\
&\cong \oplus_{n=0}^rH^0(S,{\mathcal O}_S((r-n)r_2alH+nr_2(D+alH)))\\
&\cong \sum_{m+n=r}R_{m,n}t^mu^n.
\end{array}
$$
Thus $R$, with the above grading, is the coordinate ring of $X$,
with respect to an embedding in projective space given by
$H^0(X,{\mathcal L})$.
 Since $X$ is nonsingular, the singular locus of
$\Spec(R)$ is the bigraded maximal ideal
$\sum_{m+n>0}R_{m,n}t^mu^n$ of $R$.
\end{proof}

Choose a surjective homomorphism $S_0\to R_0$ of graded ${\bf
k}$-algebras, where $S_0$ is a polynomial ring of dimension $d$.
Then the graded version of the local duality theorem
(\cite[Theorem 3.6.19]{BH}) and property (5) of $R$ imply that
\[
\Ext_{S_0}^{d-2}(R_n, S_0)\iso \left\{\begin{array}{ll}
{\bf k}^2&\text{ if $n$ is even,}\\
0&\text{ if $n$ is odd.}
\end{array}\right.
\]
More generally, if $R$ is a positively graded Noetherian
$R_0$-algebra, $M$ is a finitely generated graded $R$-module, and
$N$ an $R_0$-module, one could ask whether the graded $R$-modules
$\Ext^i_{R_0}(M,N)$, $\Ext^i_{R_0}(N,M)$ and $\Tor^i_{R_0}(M,N)$
behave tamely. The above example shows that this is in general not
the case for $\Ext^i_{R_0}(M,N)$, while for the other two functors
this is the case. In fact, computing $\Ext^i_{R_0}(N,M)$ and
$\Tor^i_{R_0}(M,N)$ by using a graded minimal free
$R_0$-resolution of $N$, these two homology groups are
subquotients of a complex whose chains are a finite number of
copies of $M$. Hence $\Ext^i_{R_0}(N,M)$ and $\Tor^i_{R_0}(M,N)$
are finitely generated graded $R$-modules. Say, $n_0$ is the
highest degree of a generator of $\Ext^i_{R_0}(N,M)$. Then one has
\[
\Ext^i_{R_0}(N,M_n)=0 \quad\text{for all $n\geq n_0$},\quad \text{
or else}\quad \Ext^i_{R_0}(N,M_n)\neq 0 \quad\text{for all $n\geq
n_0$.}
\]
The same holds true for $\Tor^i_{R_0}(M,N)$.

\section{An example of a non-tame cohomology module}

In this section we use the result of the previous section to
produce an example which gives a negative answer to the tameness
problem \cite[Problem 5.1]{BrHe} of Brodmann and Hellus. The
construction is based on a duality theorem for bigraded modules
which is given in \cite{HR}.

Let $\bf k$ be a field, $S={\bf k}[x_1,\ldots, x_r,y_1,\ldots,
y_s]$ the standard bigraded polynomial ring. In other words, we
set $\deg x_i=(1,0)$ and $\deg y_j=(0,1)$ for all $i,j$. We denote
by $P_0=(x_1,\ldots, x_r)$ the graded maximal ideal of ${\bf
k}[x_1,\ldots, x_r]$ and by $Q_0=(y_1,\ldots, y_s)$ the graded
maximal ideal of ${\bf k}[y_1,\ldots, y_s]$, and set $P=P_0S$,
$Q=Q_0S$ and $S_+=P+Q$.

Let $M$ be a finitely generated bigraded $S$-module. We set
$M_j=\Dirsum_{i\in\ZZ}M_{(i,j)}$. Then $M=\Dirsum_{j\in\ZZ}M_{j}$,
where each $M_j$ is a finitely generated graded ${\bf
k}[x_1,\ldots,x_r]$-module.

We denote by $M^\vee$ the bigraded ${\bf k}$-dual of $M$, i.e.\
the bigraded ${\bf k}$-module with components
\[
(M^\vee)_{(i,j)}=\Hom_{\bf k}(M_{(-i,-j)},{\bf k})\quad \text{for
all}\quad  i,j.
\]
By the local duality theorem one has natural isomorphisms of
bigraded $S$-modules
\[
H^i_{S_+}(M)^\vee\iso \Ext^{r+s-i}(M,S(-r,-s)) \quad \text{for
all}\quad i.
\]
In particular, all the modules $H^i_{S_+}(M)^\vee$ are finitely
generated bigraded  $S$-modules, see \cite[Theorem
3.6.19]{BH} for a similar statement in the graded case.

We shall use the following result \cite[Proposition 2.5]{HR}

\begin{Proposition}
\label{use} Suppose $M$ is a finitely generated graded generalized
Cohen-Macaulay $S$-module of dimension $d$ (i.e.\ $M$ is
Cohen-Macaulay on the punctured spectrum of $S$, or equivalently,
$H^i_{S_+}(M)$ has finite length for $i<d$). We let $N$ be the
finitely generated bigraded $S$-module $H^d_ {S_+}(M)^\vee$.  Then one
has the following long exact sequence of bigraded $S$-modules
\begin{eqnarray*}
0 \rightarrow H^{1}_P(N) \rightarrow H^{d-1}_Q(M)^\vee \rightarrow
H_{S_+}^{d-1}(M)^\vee \rightarrow H^{2}_P(N)\rightarrow
H^{d-2}_Q(M)^\vee \rightarrow H_{S_+}^{d-2}(M)^\vee \cdots
\end{eqnarray*}
\end{Proposition}

Note that $(H_{S_+}^{d-i}(M)^\vee)_j=0$ for $i>0$ and all $j\gg
0$. Thus the long exact sequence of Proposition \ref{use} yields
the following isomorphisms
\begin{eqnarray}
\label{iso} (H^{d-i}_Q(M)_{-j})^\vee\iso (H^{d-i}_Q(M)^\vee)_j\iso
H^i_P(N)_j\iso H^i_{P_0}(N_j)
\end{eqnarray}
for all $i>0$ and all $j\gg 0$.

Now let $R$ be the bigraded ${\bf k}$-algebra of Theorem
\ref{main}. We choose a bigraded presentation $S\to R$ with
$S={\bf k}[x_1,\ldots, x_r,y_1,\ldots, y_s]$, and view $R$ a
bigraded $S$-module.

We have  $\dim R=4$, so that $\omega_R=H^4_{S_+}(R)^\vee$ is the
canonical module of $R$. Since $R$ is a domain, the canonical
module localizes, that is, we have $(\omega_R)_{\wp}\iso
\omega_{R_\wp}$ for all $\wp\in \Spec(R)$, see for example
\cite[Korollar 5.25]{HK}. Furthermore,  since the singular locus
of $R$ is the bigraded maximal ideal $\mm$ of $R$, it follows that
$R_\wp$ is regular for all $\wp\neq \mm$. In particular,
$(\omega_R)_\wp\iso R_\wp$ is Cohen-Macaulay for all $\wp\neq
\mm$. This shows that $\omega_R$ is a generalized Cohen-Macaulay
module. Finally, since $R$ is normal, $R$ is in particular $S_2$,
and this, by a result of Aoyama \cite[Proposition 2]{A}, implies
that $R\iso H^4_{S_+}(\omega_R)^\vee$. Thus if we set
$J=QR=\Dirsum_{n>0}R_n$, then (\ref{iso}) applied to $\omega_R$
yields
\[
(H^2_J(\omega_R)_{-j})^\vee\iso H^2_I(R_j) \quad \text{for} \quad
j\gg 0.
\]
Thus in view of property (5) of $R$ we obtain

\begin{Corollary}
\label{counterexample} For all $j\ll 0$ one has
\[
H^2_J(\omega_R)_j\iso\left\{\begin{array}{ll}
{\bf k}^2&\text{ if $j$ is even}\\
0&\text{ if $j$ is odd.}
\end{array}\right.
\]
\end{Corollary}

Localizing at the graded maximal ideal of $R_0$ we may as well
assume that $R_0$ is local and obtain the same result.

\vskip.5truein
%\noindent
%Department of Mathematics

%\noindent University of Missouri

%\noindent Columbia, MO  65211

\end{document}